\documentclass[12pt]{article}
\usepackage{latexsym}
\usepackage{amsmath}
\usepackage{amssymb}
\usepackage{amscd}
\usepackage{graphicx}
\newtheorem{Th}{Theorem}[section]
\newtheorem{Cor}{Corollary}[section]
\newtheorem{Prop}{Proposition}[section]
\newtheorem{Lem}{Lemma}[section]
\newtheorem{Def}{Definition}[section]
\newtheorem{Rem}{Remark}[section]
\newtheorem{Ex}{Example}[section]

\newcommand{\bet}{\begin{Th}}
\newcommand{\ent}{\stepcounter{Cor}
   \stepcounter{Prop}\stepcounter{Lem}\stepcounter{Def}
   \stepcounter{Rem}\stepcounter{Ex}\end{Th}}

\newcommand{\bec}{\begin{Cor}}
\newcommand{\enc}{\stepcounter{Th}
   \stepcounter{Prop}\stepcounter{Lem}\stepcounter{Def}
   \stepcounter{Rem}\stepcounter{Ex}\end{Cor}}
\newcommand{\bep}{\begin{Prop}}
\newcommand{\enp}{\stepcounter{Th}
   \stepcounter{Cor}\stepcounter{Lem}\stepcounter{Def}
   \stepcounter{Rem}\stepcounter{Ex}\end{Prop}}
\newcommand{\bel}{\begin{Lem}}
\newcommand{\enl}{\stepcounter{Th}
   \stepcounter{Cor}\stepcounter{Prop}\stepcounter{Def}
   \stepcounter{Rem}\stepcounter{Ex}\end{Lem}}
\newcommand{\bef}{\begin{Def}}
\newcommand{\enf}{\stepcounter{Th}
   \stepcounter{Cor}\stepcounter{Prop}\stepcounter{Lem}
   \stepcounter{Rem}\stepcounter{Ex}\end{Def}}
\newcommand{\ber}{\begin{Rem} }
\newcommand{\enr}{
   \stepcounter{Th}\stepcounter{Cor}\stepcounter{Prop}
   \stepcounter{Lem}\stepcounter{Def}\stepcounter{Ex}\end{Rem}}
\newcommand{\bee}{\begin{Ex}}
\newcommand{\ene}{
   \stepcounter{Th}\stepcounter{Cor}\stepcounter{Prop}
   \stepcounter{Lem}\stepcounter{Def}\stepcounter{Rem}\end{Ex}}
\newcommand{\Proof}{\noindent{\it Proof\,}:\ }

\newcommand{\R}{\mathbf{R}}

\newcommand{\ZZ}{\mathbf{Z}}

\newcommand{\aaa}{{\boldsymbol{a}}}
\newcommand{\bbb}{{\boldsymbol{b}}}

\newcommand{\pa}{\partial}

\newcommand{\QED}{\hfill $\Box$ \par}

\newcommand{\GL}{{\mbox {\rm GL}}}

\newcommand{\Gr}{{\mbox {\rm Gr}}}

\newcommand{\ord}{{\mbox {\rm ord}}}
\newcommand{\rank}{{\mbox {\rm rank}}}

\begin{document}

\title{
Generic bifurcations of framed curves 
\\ 
in a space form and their envelopes
}
\author{Goo ISHIKAWA
\\
Department of Mathematics, Hokkaido University
}

\date{
} 

\maketitle



\section{Introduction.}

A curve endowed with a frame, called {\it a framed curve}, 
in a space-form plays important roles in topology, geometry and singularity theory. 
For example, as is well-known, the self-linking number in 
$3$-space is defined via framing (\cite{Pohl}). 
The fundamental theory of curves is formulated via osculation framing. Surface boundaries have adapted framings, etc. 
Two kinds of frames are considered in this paper; adapted frames and osculating frames 
from the viewpoint of duality. Then we classify the singularities of envelopes associated to framed curves.  
The singularities of envelopes in $E^3$ was studied in \cite{Ishikawa6} to apply to  
the flat extension problem of a surface with boundary. 
The problem on extensions by tangentially degenerate surfaces motivates to study the envelopes 
associated to framings on curves in a space form. 
In this paper we consider framed curves in 
$X = E^{n+1}$, Euclidean space, $S^{n+1}$, the sphere or $H^{n+1}$, the hyperbolic 
space of dimension $n+1$, and understand commonly in terms of projective geometry. 

Actually we work with the models  
$$
S^{n+1} = \{ x \in \R^{n+2} \mid x^2 = 1 \}, \quad 
H^{n+1} = \{ x \in \R^{1, n+1} \mid x^2 = - 1,  \ x_0 > 0 
\}, 
$$
where $\R^{1, n+1} = \R^{n+2}_1 = \{ (x_0, x_1, \dots, x_{n+1}) \}$ is the Minkowski space of 
index $(1, n+1)$ (See for instance \cite{IPS}\cite{Harvey}). The inner product in 
$\R^{1, n+1}$ is defined by $x\cdot y = - x_0y_0 + \sum_{i=1}^{n+1} x_iy_i$. 
Moreover we identify Euclidean space $E^{n+1}$ with $\{ x \in \R^{n+2} \mid x_0 = 1\} \subset \R^{n+2}$ 
if necessary. 

Let $\gamma : I \to X$ be a $C^\infty$ immersion from an interval or a circle $I$. 
In general, we mean by a {\it framing} of the immersed curve 
$\gamma$, an oriented orthonormal frame 
$(e_1, e_2, \dots, e_{n+1})$ along $\gamma$.  
%
We always 
pose the condition that 
$e_{n+1}$ is orthogonal to the velocity vector $\gamma'$. 
Then the unit normal vectors $e_{n+1}$ 
provide a $1$-parameter family of tangent hyperplanes to $\gamma$ 
and its envelope $E(\gamma)$. 

In particular, in three dimensional case ($n = 2$), 
if a framed curve $\gamma$ is given, then we have a $1$-parameter family of planes 
and its envelope surface $E(\gamma)$ in three space. 
For a $1$-parameter family of framed curves $\gamma_\lambda$, we have the $1$-parameter family of 
envelopes $E(\gamma_\lambda)$. Then we will show  

\bet
\label{envelope}
Let $\gamma_\lambda$ be a generic $1$-parameter family of framed 
curves in $E^3, S^3$ or $H^3$. 
Then the local singularity in the associated envelope $E(\gamma_\lambda)$ 
is given by one of following $5$-classes: 
{\rm (I)} the cuspidal edge, {\rm (II)} the swallowtail, 
{\rm (III)} the cuspidal breaks, {\rm (IV)} the cuspidal butterfly, 
and {\rm (V)} the full-folded-umbrella. 
\ent

In particular, the list of singularities (diffeomorphism classes), is the same for all of three geometries. 

\begin{figure}[htbp]
  \begin{center}
\includegraphics[width=5truecm, height=3truecm, clip, bb=0 520 500 850]{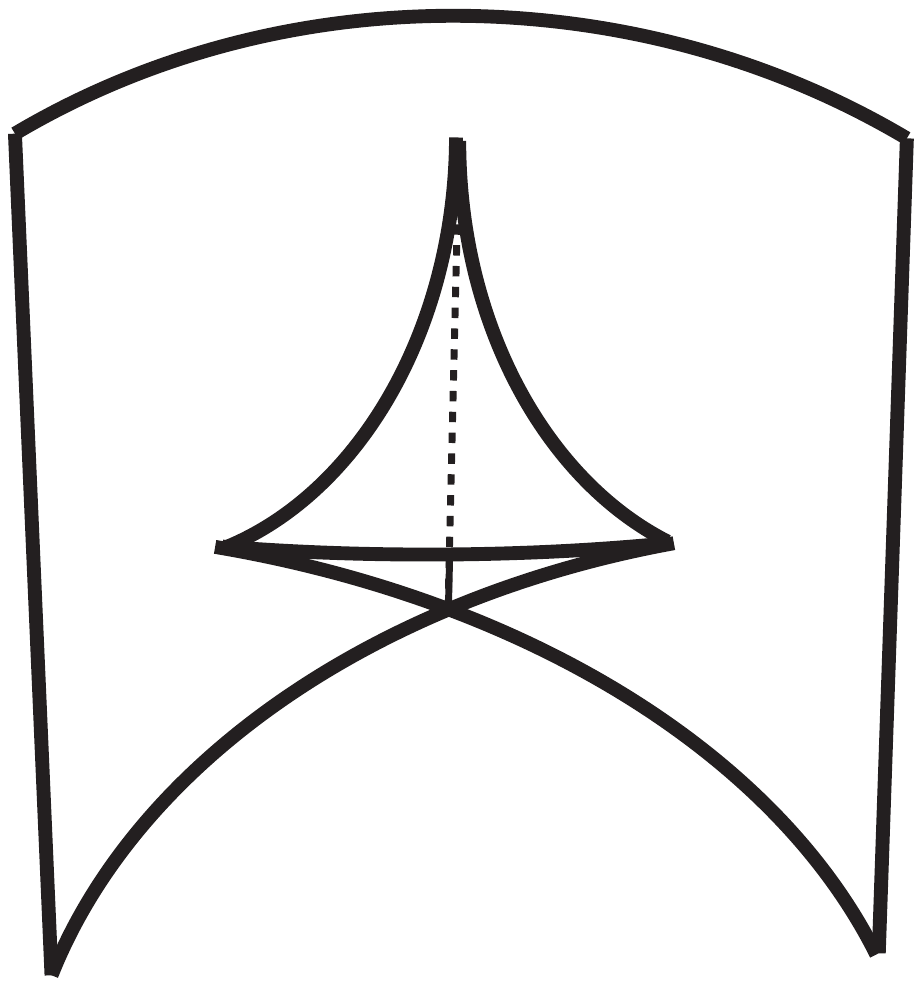}
 \hspace{2truecm}
    \includegraphics[width=5truecm, height=3truecm, clip, bb=100 440 520 800]{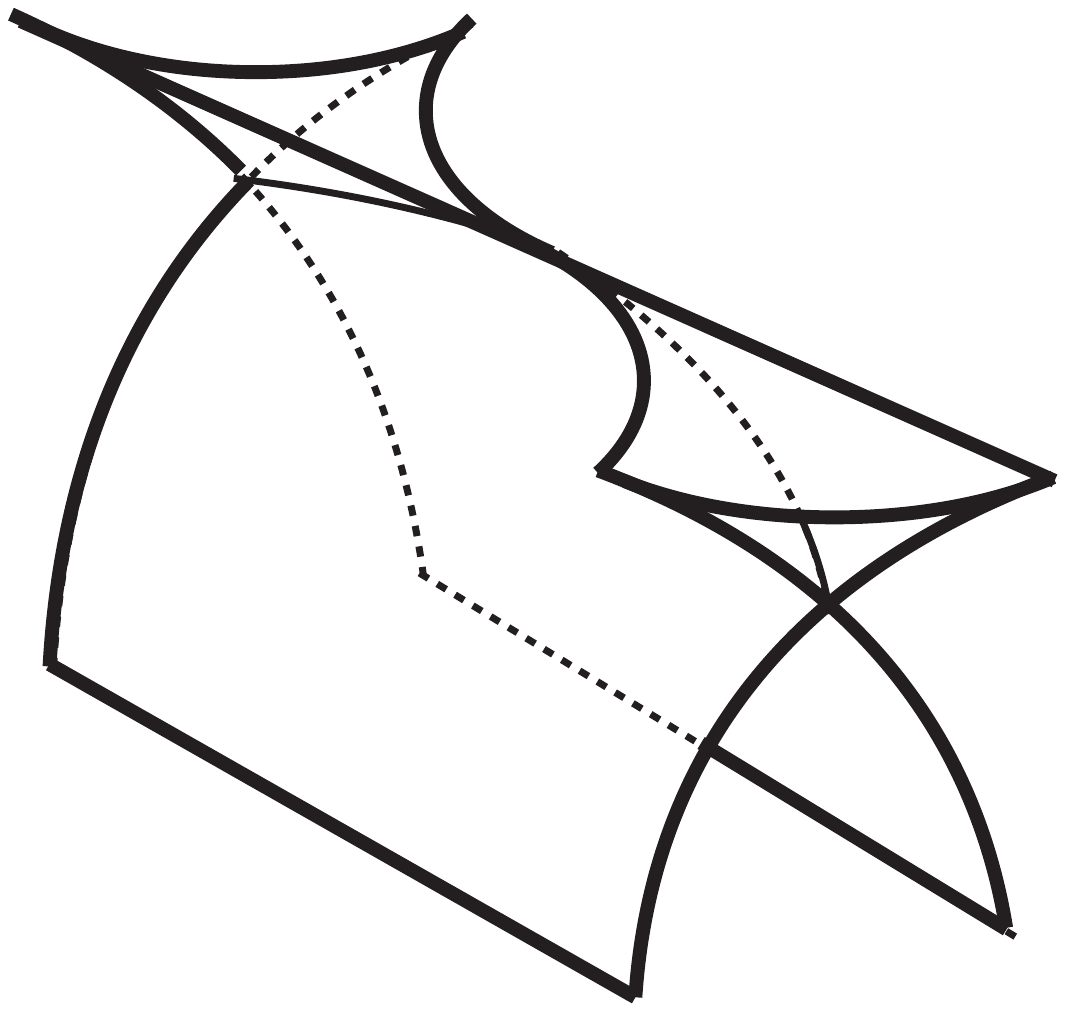}
\\
the swallowtail
and 
the cuspidal beaks
\end{center}

 \begin{center}
 \includegraphics[width=5truecm, height=3truecm, clip, bb=0 560 440 800]{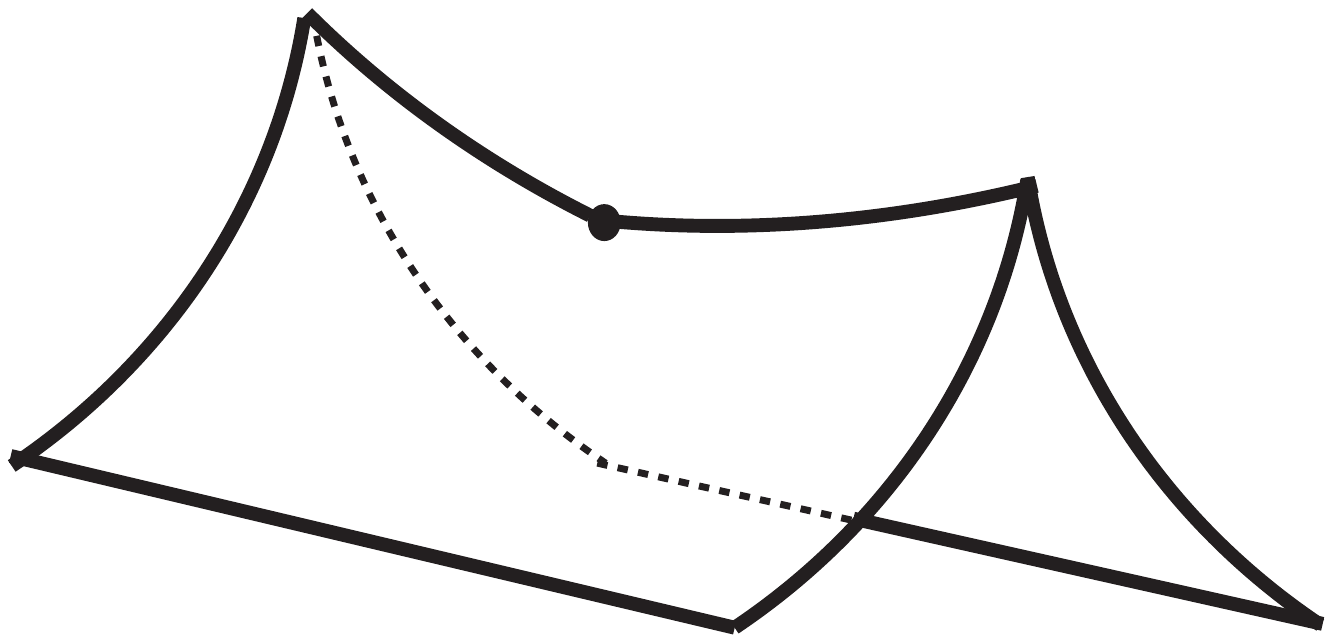}
 \hspace{2truecm}
    \includegraphics[width=5truecm, height=3truecm, clip, bb=10 400 510 800]{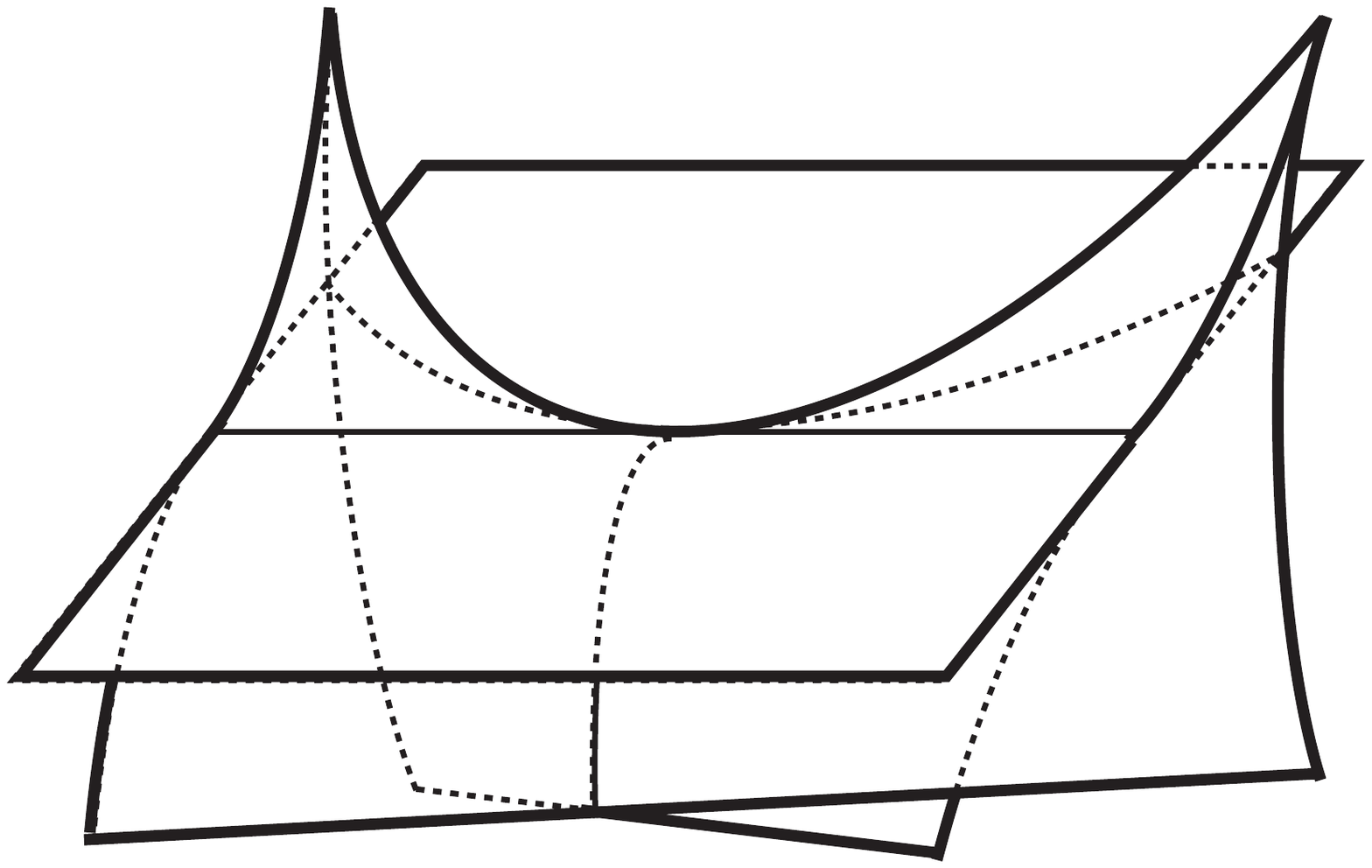}
  \\
the cuspidal butterfly 
and 
the full-folded-unbrella
\end{center}
\end{figure}%

The cuspidal edges and the swallowtails appear generically and stably. 
The swallowtails can appear in isolated positions on the envelope in any moment $\lambda$. 
The cuspidal breaks, the cuspidal butterflies, or the full-folded-umbrellas   
appear in isolated positions momentarily at isolated value $\lambda$. 
Along the parameter $\lambda$, both the cuspidal breaks and butterflies bifurcate within wavefronts (\cite{Arnold1}). However we see later they have different character in our theory. 
The full-folded-umbrella is not a wavefront (image of a non-singular Legendre submanifold), but, 
a frontal surface (image of a singular Legendre variety). 

The singularities of envelopes are closely related to singularities of tangent developables of curves. 
Tangent developables are flat in $E^3$. However they are not flat but "extrisically flat" or tangentially degenerate 
in $S^3$ and $H^3$ (cf. \cite{AG}\cite{KRUSUY}). 
In this paper the notion of {\it types} $(a_1, a_2, a_3)$ for a curve-germ is introduced and 
the cuspidal edge, (resp. the swallowtail, the cuspidal break, the cuspidal butterfly) is obtained 
as the tangent developable of a curve of type $(1, 2, 3)$ (resp. $(2, 3, 4)$, $(1, 3, 4)$, $(3, 4, 5)$). 
The cuspidal break is called also {\it Mond surface} \cite{Ishikawa4}. (See also \cite{Mond1}\cite{Mond2}). 
We have adopted the notations in \cite{IS}. 
We remark that the cuspidal butterflies bifurcate within tangent developables, however, 
the cuspidal breaks (Mond surfaces) do not. In fact we observe the Mond surface is stable 
for the deformations of curves with {\it osculating frames}. 

The full-folded-umbrella contains the tangent developable of a curve of type $(1, 2, 4)$. 
Each singularity mentioned above is given by the generating family 
$$
F(t, x_1, x_2, x_3) = 
\frac{t^{a_3}}{a_3!} + x_1\frac{t^{a_3 - a_1}}{(a_3 - a_1)!} + x_2\frac{t^{a_3 - a_2}}{(a_3 - a_2)!} + x_3 = 0
$$ 
of (totally geodesic) planes, where the normal form of the envelope is 
given by $\{ (x_1, x_2, x_3) \mid F = \frac{\pa F}{\pa t} = 0 {\mbox{\rm { for some 
}}} t \}$ (\cite{Ishikawa1}).

In this paper two kinds of frames are involved: 
one is an {\it adapted frame} of $\gamma$ which satisfies just the condition $e_1 = \gamma'$, 
the unit velocity vector field, or the differential by the arc-length parameter. 
Then $e_{n+1}$ is orthogonal to $\gamma'$. For the classification problem of envelops just that 
$e_{n+1}$ is orthogonal to $\gamma'$ is essential. 
Another is the Frenet-Serre frame of $\gamma$ along ordinary points 
where the derivatives 
$
\gamma'(t), \gamma''(t), \dots, \gamma^{(n)}(t), 
$ 
are linearly independent. Then our main idea is to introduce two kinds of distributions, or differential systems, 
on flag manifolds and regard framed curves as integral curves to those distributions. 

Bifurcations of wavefronts based on Legendre singularity theory are established by 
Arnold-Zakalyukin's theory (\cite{Arnold1}\cite{Arnold2}\cite{Zakalyukin1}\cite{Zakalyukin2}). 
The application of singularity theory to differential geometry has been developed by many authors 
(see for instance \cite{BG}). 
The singularity theory based geometry of submanifolds in hyperbolic space $H^{n+1}$ is 
initiated by Izumiya et al. (\cite{IPS}\cite{IPS2}\cite{IPT}). 
The Legendre duality developed in 
\cite{Izumiya}\cite{CI} enables us to unify the theory of framed curves in any space form as describes in this paper. 

In \S \ref{Legendre duality.}, 
we recall Legendre duality (see \cite{Bruce}\cite{IMo}\cite{CI}) 
within the level we need in this paper. 
We understand the duality in the framework of moving frames and flags in \S \ref{Moving frames and flags.}. 
After touching with non-oriented flags in \S \ref{Reduced Legendre duality.}, we introduced 
two distributions in \S \ref{Pseudo-contact and canonical distributions.}. This is very essential to study the bifurcation problem of envelopes in this paper. 
In \S \ref{Osculating flags on curves of finite type.}, the notion of type of curves is introduced 
and that of osculating flags are considered. Two kinds of framed curves are regarded as 
integral curves to two kinds of distributions. Then 
we prove codimension formulae for framed curves in \S \ref{Integral curves and codimension formula.}, 
which implies the classification results of singularities of envelopes including Theorem\ref{envelope}, 
in \S \ref{Singularities of envelopes and their bifurcations.}. 


\section{Legendre duality.} 
\label{Legendre duality.} 

Though in this paper we mainly treat curves in Riemannian spaces $X = E^{n+1}, S^{n+1}, H^{n+1}$, 
regarding the duality, naturally we work in other spaces as well. 
In particular we are led to consider {\it de Sitter space}
$$
S^{1, n} = \{ x \in \R^{1, n+1} \mid x^2 = 1\}, 
$$
which is a semi-Riemannian manifold, 
since any vector 
of an frame $(e_1, \dots, e_{n+1})$ along 
a curve in $H^{n+1}$ belongs to $S^{1, n}$. 

We regard $\widehat{\gamma} = e_{n+1}$ a curve in $Y = S^{n+1}$ 
(resp. in $Y = S^{1, n}$) 
if $X = S^{n+1}$ (resp. $X = H^{n+1}$). 
In Euclidean case, we set 
$\widehat{\gamma} = (- \gamma\cdot e_{n+1}, e_{n+1})$ 
and regard it as a curve in $Y = \R\times S^{n+1}$. 
We call $\widehat{\gamma}$ the {\it frame dual} to $\gamma$. 
Then, in any case, the \lq\lq type" of the curve $\widehat{\gamma}$ in $Y$ 
describes the singularities of the envelope 
$E(\gamma)$ in $X$.  

We denote by $Z = {\widetilde{\Gr}}(n, TX)$ the manifold of oriented tangent hyperplanes of $X$ 
and by $\pi_1 : Z \to X$ the projection which maps a hyperplane $\Pi \subset T_xX$ to $x \in X$. 
A framed curve $\gamma : I \to X$ with the framing $(e_1, \dots, e_{n+1})$ 
lifts to a curve $\widetilde{\gamma} : I \to Z$ which is 
defined by $\widetilde{\gamma}(t) = \langle e_1(t), \dots, e_n(t) \rangle_{\R}$. 
In each of three cases, $Z$ is identified with the unit tangent bundle $T_1X$ via the metric, 
actually with 
$T_1E^{n+1} = E^{n+1}\times S^n$, 
$$
\begin{array}{rcl}
T_1S^{n+1} & = & \{ (x, y) \in S^{n+1}\times S^{n+1} \mid x\cdot y = 0 \}, \quad {\mbox{\rm and}}, \\
T_1H^{n+1} & = & \{ (x, y) \in H^{n+1}\times S^{1, n} \mid x\cdot y = 0 \}.  
\end{array}
$$
Then, under the above identification, 
the lifting $\widetilde{\gamma} : I \to Z$ is given by $\widetilde{\gamma}(t) = e_{n+1}(t)$ (\cite{IM2}). 

Consider the contact structure on $Z$: 
the one-form $\theta = v\cdot dx$ restricted to $Z = E^{n+1}\times S^n$, 
$\theta = y\cdot dx$ restricted to $Z = T_1S^{n+1}$ or $T_1H^{n+1}$, 
is a contact form on $Z$. 
In elliptic or hyperbolic case, let $\pi_2 : Z \to Y$ be the projection defined by $\pi_2(x, y) = y$.  
In Euclidean case, let $\pi_2 : Z \to Y$ be the projection 
defined by $\pi_2(x, y) = (-x\cdot y, \ y) (x \in E^{n+1}, y \in S^n)$. 
Then we see both $\pi_1 : Z \to X$ and $\pi_2 : Z \to Y$ are Legendre fibrations. 

Suppose the framing of $\gamma : I \to X$ satisfies the condition $e_1 = \gamma'$. 
Then $e_{n+1}$ is normal to $\gamma'$. 
And then we see that the lifting 
$\widetilde{\gamma} : I \to Z$ of $\gamma : I \to X$ turns to be 
{\it integral} in the sense that ${\widetilde{\gamma}}^* \theta = 0$. 
The lifting $\widetilde{\gamma} : I \to Z$ of a framed immersion   
$\gamma : I \to X$ defines 
a \lq\lq sub-front" $\widehat{\gamma} = 
\pi_2\circ\widetilde{\gamma} : I \to Y$ possibly with singularities, 
in the sense that the integral lifting $\widetilde{\gamma}$ with respect to $\pi_2$ 
is attached to the just parametrised curve $\widehat{\gamma}$ in $Y$. 

Note that, in the case $X = H^{n+1}$, $Z = T_1 H^{n+1}$ is identified with 
$T_{-1} S^{1, n}$, the manifold of tangent vectors $v \in T_y S^{1, n}$ with $v^2 = -1$ 
to the semi-Riemannian manifold $S^{1, n}$ (\cite{IM2}). 


As the model of duality, we do have the projective duality (\cite{O.P.Scherbak1}\cite{IMo}); we set 
$$
Z \ = \ {\mathcal I}_{n+2} \ := \ \{ ([x], [y]) \in P^{n+1}\times P^{n+1*} \mid x\cdot y = 0 \}. 
$$
Here $P^{n+1*}$ is the dual projective space and 
$\cdot$ means the natural paring. The contact structure on ${\mathcal I}_{n+2}$ 
is defined by $dx\cdot y = x\cdot dy = 0$ (\cite{IMo}). 
The projections $\pi_1 : {\mathcal I}_{n+2} \to X = P^{n+1}, \pi_2 : {\mathcal I}_{n+2} \to Y = P^{n+1*}$ 
are both Legendre fibrations. 

The following fact is basic to unify our treatment: 

\bep
\label{flatness}
{\rm (\cite{IM}\cite{IM2})}
All Legendre double fibrations $X \longleftarrow Z \longrightarrow Y$ constructed above 
are locally isomorphic to each other. In particular each of them is locally isomorphic to 
the double fibration of the projective duality $P^{n+1} \longleftarrow {\mathcal I}_{n+2} \longrightarrow P^{n+1*}$. 
\enp

The proof of Proposition \ref{flatness} is proved naturally via the underlying 
flag structure that we are going to explain. 



\section{Moving frames and flags.}
\label{Moving frames and flags.}

For a framed curve $\gamma : I \to X$ with a frame 
$(e_1, \dots, e_{n+1})$, naturally there is associated a 
\lq\lq moving frame" $\widetilde{\gamma} : I \to G$, for each case, 
in a Lie subgroup of $\GL_+(n+2, \R)$, regular matrices with positive determinant. 

If $X = E^{n+1}$, then we set $e_0(t) = \gamma(t) \in E^{n+1}$, 
and we have the moving frame 
$\widetilde{\gamma} = (e_0, e_1, \dots, e_{n+1}) : I \to G = {\mbox{\rm{Euc}}}(E^{n+1}) 
\subset \GL_+(n+2, \R)$ in 
the group of orientation preserving Euclidean motion on $E^{n+1}$. 
If $X = S^{n+1}$, then we set $e_0(t) = \gamma(t) \in S^{n+1}$, 
and we have the moving frame 
$\widetilde{\gamma} = (e_0, e_1, \dots, e_{n+1}) : I \to G = SO(n+2) \subset \GL_+(n+2, \R)$. 
If $X = H^{n+1}$, then we set $e_0(t) = \gamma(t) \in H^{n+1}$, 
and we have the moving frame 
$\widetilde{\gamma} = (e_0, e_1, \dots, e_{n+1}) : I \to G = SO(1,n+1) \subset \GL_+(n+2, \R)$. 
In any of three cases, the frame manifold $G$ is identified with an open subset of 
the oriented flag manifold $\widetilde{\mathcal F}_{n+2}$ 
consisting of oriented complete flags 
$$
V_1 \subset  V_2  \subset \cdots \subset V_{n+1} \subset \R^{n+2} 
$$
in $\R^{n+2}$. For each $g = (e_0, e_1, \dots, e_{n+1}) \in \GL_+(n+2, \R)$, 
we set the oriented subspace 
$$
V_i = \langle e_0, e_1, \dots, e_{i-1}\rangle_{\R} \subset \R^{n+2}, \ (1 \leq i \leq n+1). 
$$
This induces an open embedding $G \to \widetilde{\mathcal F}_{n+2}$. 
Note that $\widetilde{\mathcal F}_{n+2}$ is the quotient  
space of $\GL_+(n+2, \R)$ by upper triangular matrices, and 
$\dim G = \dim \widetilde{\mathcal F}_{n+2} = \frac{(n+1)(n+2)}{2}$. 
Moreover note that the inner product restricted to each $V_i$ is non-degenerate. 
Therefore $G$ is embedding in non-degenarate flags $\widetilde{\mathcal F}_{n+2}^{0} \subset \widetilde{\mathcal F}_{n+2}$ consisting of flags $(V_1, \dots, V_{n+1})$ where 
the inner product restricted to each $V_i$ is non-degenerate. Remark that 
$\widetilde{\mathcal F}_{n+2}^{0}$ is open dense in $\widetilde{\mathcal F}_{n+2}$. 
However note that in \cite{Izumiya}\cite{CI}\cite{IS}, more general framings are considered to treat also 
{\it the light cone} in Minkowski space.




Thus, for a framed curve $\gamma : I \to X$ in $X = E^{n+1}, S^{n+1}, H^{n+1}$, 
with the frame $(e_1, \dots, e_{n+1})$, 
we have the \lq\lq framed curve" $\widetilde{\gamma} : I \to \widetilde{\mathcal F}_{n+2}^{0}$ 
by setting 
$$
V_i(t)  = \langle e_0(t), e_1(t), \dots, e_{i-1}(t) \rangle_{\R} \subset \R^{n+2}, \ (1 \leq i \leq n+1). 
$$
Then $\widetilde{\gamma}$ is a lifting of 
$\gamma$ for the projection $\pi_1 : \widetilde{\mathcal F}_{n+2}^{0} \to \widetilde{\Gr}(1, \R^{n+2})$
to Grassmannian of oriented lines in $\R^{n+2}$. Note that 
there is the natural open embedding $X \subset \widetilde{\Gr}(1, \R^{n+2})$ in each of three cases. 

\section{Reduced Legendre duality.}
\label{Reduced Legendre duality.}

The vector space $\R^{n+2}$ have $\ZZ/2\ZZ$-action defined by 
$x \mapsto - x$. To describe the duality we are going to use properly, 
it is natural to take quotient and set   
$$
\begin{array}{cc}
AG(n, n+1) \ := \ \{ (r, y) \in \R\times S^n \}/(\ZZ/2\ZZ), \quad 
P^{n+1} \ := \ \{ x \in \R^{n+2} \mid x^2 = 1\}/(\ZZ/2\ZZ), 
\vspace{0.2truecm}
\\
H^{n+1} \ := \ \{ x \in \R^{1, n+1} \mid x^2 = -1\}/(\ZZ/2\ZZ), 
{\mbox{\rm 
\ and, \ }} 
P^{1, n} \ := \ \{ x \in \R^{1, n+1} \mid x^2 = 1\}/(\ZZ/2\ZZ). 
\end{array}
$$
We call $P^{n+1}$ the {\it elliptic space} and $P^{1, n}$ the {\it reduced de Sitter space}. 
Remark that $AG(n, n+1)$ is identified with the set of affine non-oriented hyperplanes in $E^{n+1}$. 
We regard $P^{n+1}$ (resp. $P^{1, n}$) as the double-quotient of the sphere $S^{n+1}$ (resp. 
$S^{1, n}$) with the induced metric. 
Set 
$X = E^{n+1}, P^{n+1}, H^{n+1}$ in Euclidean, elliptic, hyperbolic case,  
respectively. Then we set $Y = AG(n, n+1), P^{n+1}, P^{1, n}$ respectively. 

We consider the incidence manifold in each geometry: 
$$
Z \ := \ \{ ([x], [y]) \in X\times Y \mid x\cdot y = 0 \}, 
$$
for elliptic and hyperbolic cases, 
and 
$$
Z \ := \ \{ (x, [r, y]) \in X\times Y \mid x\cdot y + r = 0 \}, 
$$
for Euclidean case. 
In each case, $Z$ is regarded naturally as an open subset of $PT^*X$ 
and is endowed with the standard contact structure. 
Then the double fibrations 
$\pi_1 : Z \to X$ and $\pi_2 : Z \to Y$ are Legendre. Moreover all 
Legendre double fibrations $X \longleftarrow Z \longrightarrow Y$ are 
locally isomorphic 
to the projective duality $P^{n+1} \longleftarrow {\mathcal I}_{n+2} \longrightarrow P^{n+1*}$ as in 
Proposition \ref{flatness}. 

Let ${\mathcal F}_{n+2}$ be the manifold of non-oriented complete flags 
$$
V_1 \subset  V_2  \subset \cdots \subset V_{n+1} \subset \R^{n+2}, 
$$
consisting of vector subspaces $V_i$ of dimension $i$ in $\R^{n+2}$. 
The forgetful mapping $\pi : 
\widetilde{\mathcal F}_{n+2} \to {\mathcal F}_{n+2}$ forms a covering of order 
$2^{n+1}$. 
For a framed curve in a reduced space, the lifting is a curve in a non-oriented flag manifold. 


\section{Pseudo-contact and canonical distributions.} 
\label{Pseudo-contact and canonical distributions.}

We  will consider two classes of curves in the frame manifold $\widetilde{\mathcal F}_{n+2}$, 
by introducing two kinds of distributions 
$\widetilde{\mathcal C} \subset \widetilde{\mathcal D} \subset T\widetilde{\mathcal F}_{n+2}$. 
Denote by $\pi_i : \widetilde{\mathcal F}_{n+2} \to \widetilde{\Gr}(i, \R^{n+2})$ the projection to 
Grassmannian of oriented $i$-planes in $\R^{n+2}$ defined by 
$
\pi_i(V_1, \dots, V_i, \dots, V_{n+1}) = V_i. 
$
Then we define, 
for $v \in T\widetilde{\mathcal F}_{n+2}$, 
$v \in \widetilde{\mathcal D}_{(V_1, \dots, V_{n+1})}$ if 
${\pi_1}_*(v) \in T\widetilde{\Gr}(1, V_{n+1}) (\subset T\widetilde{\Gr}(1, \R^{n+2}))$, while 
$v \in \widetilde{\mathcal C}_{(V_1, \dots, V_{n+1})}$ if 
${\pi_i}_*(v) \in T\widetilde{\Gr}(i, V_{i+1}) (\subset T\widetilde{\Gr}(i, \R^{n+2})), (1 \leq i \leq n)$. 

We call the distribution $\widetilde{\mathcal D}$ {\it pseudo-contact distribution} and 
$\widetilde{\mathcal C}$ {\it canonical distribution}. 
Note that 
the rank of $\widetilde{\mathcal C}$ (resp. 
$\widetilde{\mathcal D}$) is $n+1$ (resp. $\frac{(n+1)(n+2)}{2} - 1$) in 
$T\widetilde{\mathcal F}_{n+2}$. 
Both $\widetilde{C}$ and $\widetilde{D}$ are bracket generating; in fact, $n$-th bracket $\widetilde{\mathcal C}^n$ 
of $\widetilde{\mathcal C}$ coincides with $\widetilde{\mathcal D}$.  
Denote by $\widetilde{\mathcal I}_{n+2}$ the flag manifold 
consisting of flag $V_1 \subset V_{n+1} \subset \R^{n+2}$ with an oriented line $V_1$ and 
an oriented hyperplanes $V_{n+1}$. 
Consider the canonical projection 
$
\pi_{1, n+1} : \widetilde{\mathcal F}_{n+2} \longrightarrow \widetilde{\mathcal I}_{n+2}
$
defined by $\pi_{1, n+1}(V_1, V_2, \dots, V_{n+1}) = (V_1, V_{n+1})$. 
Then we have $\widetilde{\mathcal D} = (\pi_{1, n+1})_*^{-1}(D)$, 
the pull-back of the contact structure $D$ on 
$\widetilde{\mathcal I}_{n+2}$:  
for $v \in T\widetilde{\mathcal I}_{n+2}$, $v \in \widetilde{D}_{(V_1, V_{n+1})}$ if $(\pi_1)_*(v) 
\in T\widetilde{\Gr}(1, V_{n+1})$. The contact structure $D$ on $\widetilde{\mathcal I}_{n+2}$ 
is the pull-back of the contact structure on ${\mathcal I}_{n+2}$ introduced in \S \ref{Legendre duality.}. 

Similar constructions go as well for non-oriented case. 

Define two distributions (vector sub-bundles) ${\mathcal C} \subset {\mathcal D} 
\subset T{\mathcal F}_{n+2}$ 
on the non-oriented flag manifold ${\mathcal F}_{n+2}$ 
as follows: For $v \in T{\mathcal F}_{n+2}$, 
$v \in {\mathcal D}_{(V_1, \dots, V_{n+1})}$ if 
${\pi_1}_*(v) \in T{\Gr}(1, V_{n+1}) (\subset T{\Gr}(1, \R^{n+2}))$, while 
$v \in {\mathcal C}_{(V_1, \dots, V_{n+1})}$ if 
${\pi_i}_*(v) \in T{\Gr}(i, V_{i+1}) (\subset T{\Gr}(i, \R^{n+2})), (1 \leq i \leq n)$. 
We call also the distribution ${\mathcal D}$ {\it pseudo-contact distribution} and 
${\mathcal C}$ {\it canonical distribution}. 
Clearly the forgetful covering $\pi : \widetilde{\mathcal F}_{n+2} \to {\mathcal F}_{n+2}$ 
induces a local isomorphism of $\widetilde{\mathcal D}$ and ${\mathcal D}$ (resp. 
$\widetilde{\mathcal C}$ and ${\mathcal C}$). The pseudo-contact structure ${\mathcal D}$ 
is the pull-back of the contact structure on ${\mathcal I}_{n+2}$ 
via the canonical projection $\pi_{1, n+1} : {\mathcal F}_{n+2} \to {\mathcal I}_{n+2}$. 

Now we describe the local structure of the canonical distribution 
${\mathcal C} \subset T{\mathcal F}_{n+2}$. 
Since ${\mathcal C}$ is $\GL(n+2, \R)$-invariant, we describe ${\mathcal C}$ 
in a neighbourhood of the standard flag $E \in {\mathcal F}_{n+2}$ 
which corresponds to the unit matrix. 
The flag manifold has local coordinates $x_i^{\ j}, (0 \leq j < i \leq n+1)$ near $E$ as components of lower 
triangular matrices. 
Then ${\mathcal C}$ is defined by the system of $1$-forms 
$$
dx_i^{\ j}  -  x_i^{\ j+1} dx_{j+1}^{\ j} = 0, \qquad (0 \leq j, \ j+1 < i). 
$$
Therefore a ${\mathcal C}$-integral curve $\Gamma(t) = (x_i^{\ j}(t))_{0 \leq j < i \leq n+1}$ 
through the standard flag $E \in {\mathcal F}_{n+2}$ 
is determined just by $x_j^{j-1}(t), (1 \leq j \leq n+1)$.

\ber
{\rm
The complete flag manifold ${\mathcal F}_{n+2} = {\mathcal F}(\R^{n+2})$ (resp. 
$\widetilde{\mathcal F}_{n+2} = \widetilde{\mathcal F}(\R^{n+2})$) possesses the duality 
between 
${\mathcal F}^*_{n+2} = {\mathcal F}(\R^{n+2*})$ (resp. 
$\widetilde{\mathcal F}^*_{n+2} = \widetilde{\mathcal F}(\R^{n+2*})$)
by 
$$
(V_1, V_2, \dots, V_n, V_{n+1}) \mapsto 
(V_{n+1}^{\vee}, V_n^{\vee}, \dots, V_2^{\vee}, V_1^{\vee})
$$
where $V^{\vee} \subset \R^{n+2*}$ is the annihilator for $V \subset \R^{n+2}$. 
Then, for each metric on $\R^{n+2}$, the dual space $\R^{n+2*}$ is identified with 
$\R^{n+2}$. Thus we have the canonical involution on ${\mathcal F}(\R^{n+2})$ (resp. 
$\widetilde{\mathcal F}(\R^{n+2}), {\mathcal F}(\R^{1, n+1}), \widetilde{\mathcal F}(\R^{1, n+1})$). 
Similarly we have the canonical involution on 
${\mathcal I}(\R^{n+2})$ (resp. 
$\widetilde{\mathcal I}(\R^{n+2}), {\mathcal I}(\R^{1, n+1}), \widetilde{\mathcal I}(\R^{1, n+1})$). 
This justifies our theory. 
}
\enr


\section{Osculating flags on curves of finite type.}
\label{Osculating flags on curves of finite type.}

In general we treat a curve {\it of finite type} and define an analogue of Frenet-Serre frame 
even when the curve is not an immersion. 
Here, since $X \subset \R^{n+2} \setminus \{ 0\}$, we regard $\gamma$ as a curve in 
$\R^{n+2} \setminus \{ 0\}$.  The metric on $\R^{n+2}$ does not concern here. 

Let $\gamma : I \to \R^{n+2} \setminus \{ 0\}$ be a $C^\infty$ curve. 
The curve $\gamma$ is called of finite type at $t = t_0 \in I$ 
if the $(n+2)\times \infty$-matrix 
$$
\tilde{A}(t) = 
\left(
\gamma(t), \gamma'(t), \gamma''(t), \cdots, \gamma^{(r)}(t), \cdots 
\right), 
$$
is of rank $n+2$ for $t = t_0$. 
We set $(n+2)\times(r+1)$-matrix 
$$
\tilde{A}_r(t) = 
\left(
\gamma(t), \gamma'(t), \gamma''(t), \cdots, \gamma^{(r)}(t)
\right). 
$$
Then $\gamma$ is of finite type at $t = t_0$ if $\tilde{A}_r(t)$ is of rank $n+2$ for a sufficiently 
large $r$.  

Let $\aaa = (a_1, \dots, a_n, a_{n+1})$ be a sequence of strictly increasing natural numbers, 
$1 \leq a_1 < \cdots < a_n < a_{n+1}$. 
Then we call $\gamma$ of type $\aaa$ at $t = t_0 \in I$ if 
$$
\begin{array}{c}
\min \{ r \mid \rank \tilde{A}_r(t_0) = 2\} = a_1,  \quad 
\min \{ r \mid \rank \tilde{A}_r(t_0) = 3\} = a_2, \quad \dots,  
\vspace{0.2truecm}
\\
\min \{ r \mid \rank \tilde{A}_r(t_0) = n+2\} = a_{n+1}. 
\end{array}
$$

We can define type for curves in the reduced space $P^{n+1}$ as well, 
by just considering the double covering. 

A point $\gamma(t_0)$ on $\gamma$ is called an {\it ordinary point} if 
$\gamma$ is of type $(1, 2, \dots, n+1)$. 
Otherwise it is called a {\it special point}. 
The parameters of special points form discrete subset in $I$ if $\gamma$ is of finite type. 

If $\gamma$ is of type $\aaa$ at $t = t_0$, then we set 
$$
O_i(t_0) = \langle \gamma(t_0), \gamma'(t_0), \dots, \gamma^{(a_{i-1})}(t_0) \rangle_{\R}, 
$$
which is, by definition, 
an $i$-dimensional subspace of $\R^{n+2}$,  $(1 \leq i \leq n+1)$. 
Then we have 
\bel
\label{Osculating flag}
The curve $\widetilde{\gamma} : I \to {\mathcal F}_{n+2}$ in the non-oriented 
flag manifold ${\mathcal F}_{n+2}$ defined by 
$$
\widetilde{\gamma}(t) : O_1(t) \subset O_2(t) \subset \dots \subset O_{n+1}(t) \subset \R^{n+2}
$$ 
is a $C^\infty$ curve. 
Moreover we can give an orientation on the flag 
locally near $t_0 \in I$. Namely we have local lifting of $\gamma$ in $\widetilde{\mathcal F}_{n+2}$ 
for the forgetful covering  
$\pi: \widetilde{\mathcal F}_{n+2} \to {\mathcal F}_{n+2}$ from the manifold of oriented flags to those of non-oriented flags. 
\enl

We call $\widetilde{\gamma}(t)$ the {\it osculating flag} of $\gamma$ at $t \in I$.  

\

\noindent{\it Proof of Lemma \ref{Osculating flag}.} 
Consider $(n+2)\times(n+2)$-matrix $B(t) = (\gamma(t), \gamma^{(a_1)}, \dots, \gamma^{(a_{n+1})})$. 
We may suppose, after a suitable linear transformation of $\gamma$ in $\R^{n+2}$, that 
$B(t_0)$ is the unit matrix. Then the lower triangular components of the matrix 
$(\gamma(t), \gamma'(t), \dots, \gamma^{(n+1)}(t))$ provides the local representation of 
$\widetilde{\gamma}$ 
in terms of local coordinates $\widetilde{\mathcal F}_{n+2}$ near $\widetilde{\gamma}(t_0)$. 
\QED

\

Suppose $\gamma$ is a curve in $X = E^{n+1}, S^{n+1} (\subset \R^{n+2})$ or $X = H^{n+1} (\subset \R^{1, n+1})$ and moreover suppose, in the case $X = H^{n+1}$, 
the restriction of the metric to each $O_i(t)$ is non-degenerate. 
If an oriented flag field along $\gamma$ is given, 
then an orthonormal frame $(e_1, \dots, e_{n+1})$ along 
$\gamma$ is uniquely constructed by the Gram-Schmidt's orthogonalisation 
which depends on the given orientation. 
We call this frame on $\gamma$ an {\it osculating frame}. 
For instance, there exists the unique unit vector $e_1(t) \in O_2(t)$ normal to $e_0(t) = \gamma(t)$ 
such that $(e_0(t), e_1(t))$ forms an oriented basis of $O_2(t)$. 
We see, by Lemma \ref{Osculating flag}, 
any osculating frame constructed above is $C^\infty$ along $\gamma$ which coincides, up to sign pointwise, with Frenet-Serre frame on ordinary points.


%


\section{Integral curves and codimension formula.}
\label{Integral curves and codimension formula.}

For an adapted framing, the lifting $\widetilde{\gamma} : I \to G$ is an integral curve 
to $\widetilde{\mathcal D}$. Moreover, for an osculating framing, 
$\widetilde{\gamma} : I \to G$ is an integral curve to $\widetilde{\mathcal C}$. 
(See \S \ref{Pseudo-contact and canonical distributions.}.) 
Thus we regard the class of adapted framed curves as the class of $\widetilde{\mathcal D}$-integral curves 
in $G$ or $\widetilde{\mathcal F}_{n+2}$ or, as being locally equivalent, 
the class of ${\mathcal D}$-integral curves in ${\mathcal F}_{n+2}$. 

On the other hand, a curve of finite type $\gamma : I \to X$ lifts, via the osculating flag, 
to a ${\mathcal C}$-integral curve 
$\widetilde{\gamma} : I \to {\mathcal F}_{n+2}$ globally. 
Moreover $\gamma$ lifts locally to a $\widetilde{\mathcal C}$-integral curve 
$\widetilde{\gamma} : I \to \widetilde{\mathcal F}_{n+2}$, which satisfies 
$e_1 = \pm\gamma'$ (arc-length differential) on immersive points pointwise. 

\ber
{\rm 
A $C^\infty$ family of curves of finite types $\gamma_\lambda : I \to X$ 
needs not to be liftable, even locally, as a 
$C^\infty$ family of ${\mathcal C}$-integral curves 
$\widetilde{\gamma}_\lambda : I \to {\mathcal F}_{n+2}$. 
The osculating flags do not behave smoothly under arbitrary deformation of curves. 
This is why we consider also the class of ${\mathcal C}$-integral integral curves 
for the bifurcation problem of envelopes. 
}
\enr

Now we consider three kinds of jet spaces of curves. 

First, we recall the ordinary jet space $J^r(I, X)$ consisting of $r$-jets of curves $I \to X$ or 
$J^r(I, Y)$ for curves $I \to Y$. Their local description are the same as in the case $X = Y = P^{n+1}$. 
O.P.Scherbak \cite{O.P.Scherbak1} shows that the codimension, called the {\it jet-codimensionin} 
${\mbox{\rm{Jet-codim}}}(\aaa)$, in the jet space $J^r(I, P^{n+1})$ of the set $\Sigma(\aaa)$ of curves in $P^{n+1}$ 
of type 
$
\aaa = (a_1, a_2, \dots, a_{n+1})
$ 
is given, for sufficiently large $r$, by 
$$
{\mbox{\rm{Jet-codim}}}(\aaa) \ = \ s(\aaa) \ := \ \sum_{i=1}^{n+1} (a_i - i), 
$$
the Schubert number which appears in Schubert calculus (\cite{MS}\cite{Kazarian2}). 

Second, we consider the jet space of ${\mathcal D}$-integral curves, 
$J^r_{\mathcal D}(I,  {\mathcal F}_{n+2}) \subset J^r(I, {\mathcal F}_{n+2})$. 
Each ${\mathcal D}$-integral curve 
$\Gamma : I \to {\mathcal F}_{n+2}$ projects to a curve $\pi_1\circ \Gamma : 
I \to P^{n+1} = \Gr(1, \R^{n+2})$ by the canonical projection 
$\pi_1 : {\mathcal F}_{n+2} \to  P^{n+1}$, $\pi_1(V_1, V_{n+1}) = V_1$. 
Then, given type $\aaa$, we have the set of jets $\Sigma_{\mathcal D}({\aaa})$ in 
$J^r_{\mathcal D}(I,  {\mathcal F}_{n+2})$. We denote its codimension by 
${\mbox{\rm{Jet-codim}}}_{\mathcal D}(\aaa)$. 

Then we have 

\bet
\label{codim-D}
The jet-codimension of the set of ${\mathcal D}$-integral curves 
$\Gamma : I 
\to {\mathcal F}_{n+2}$ such that 
$\pi_1\circ\Gamma$ is of type $\aaa = (a_1, a_2, \dots, a_{n+1})$, is given by 
$$
{\mbox{\rm{Jet-codim}}}_{\mathcal D}(\aaa)  =  {\displaystyle \sum_{i=2}^{n+1}} (a_i - i) \ = \ s(\aaa) - (a_1 - 1).
$$
\ent

\Proof 
We may suppose the case point $t_0 = 0$. 
Consider the integral jet space $J^r_{\rm{int}}(1, 2n+1)$ on germs of  
integral curves $\Gamma : I \to 
{\mathcal I}_{n+2}$ to the contact structure. 
Denote by $\Sigma_{a_1} \subset J^r_{\rm{int}}(1, 2n+1)$ the set of integral jet $j^r\Gamma$ with 
$\pi_1\circ \Gamma$ is of order $\geq a_1$. 
Take Darboux coordinates 
$x_1, \dots, x_n, z, p_1, \dots, p_n$ of ${\mathcal I}_{n+2}$ centred at  $\Gamma(t_0)$ and 
so that the contact structure is given by $dz - (p_1 dx_1 + \cdots + p_n dx_n) = 0$ and 
$\pi_1 : {\mathcal I}_{n+2} \to {\Gr}(1, \R^{n+2})$ 
is given by $(x_1, \dots, x_n, z, p_1, \dots, p_n) \mapsto (x_1, \dots, x_n, z)$. 
Let $\Gamma(t) = (x_1(t), \dots, x_n(t), z(t), p_1(t), \dots, p_n(t))$ and 
without loss of generality, we suppose $\ord x_1 = a_1$. 
Consider the mapping $\Pi : J^r_{\rm{int}}(1, 2n+1) \to J^{r - a_1 +1}(1, n+1)$ 
defined by 
$$
\Pi(x_1, \dots, x_n, z, p_1, \dots, p_n) = (x_1/t^{a_1-1}, \dots, 
x_n/t^{a_1-1}, z/t^{a_1-1}). 
$$
Take any deformation $c(t, s) = (X_1(t, s), \dots, X_n(t, s), Z(t, s))$ of $\Pi(j^r \Gamma)$ at $s = 0$. 
We set 
$$
\begin{array}{ccl}
P_1(t, s) & := & \dfrac{(t^{a_1-1}Z(t, s))'}{(t^{a_1-1}X_1(t, s))'} - 
\sum_{i=2}^n p_i(t)\dfrac{(t^{a_1-1}X_i(t, s))'}{(t^{a_1-1}X_1(t, s))'}, 
\vspace{0.2truecm}
\\
P_i(t, s) & := & p_i(t), (i = 2, \dots, n), 
\end{array}
$$
for representatives at $(t, s) = (0, 0)$. 
Then we get the integral deformation 
$$
C(t, s) = (X_1(t, s), \dots, X_n(t, s), Z(t, s), P_1(t, s), \dots, P_n(t, s))
$$ 
of $\Gamma(t)$ at $s = 0$, which satisfies 
$\pi(C(t, s)) = c(t, s)$. This show that any curve starting at $j^k(\pi\circ \Gamma)(0)$ 
in $J^{r - a_1}(1, n+1)$ lifts to a curve starting at $j^r\Gamma(0)$ 
in $J_{\rm{int}}^k(1, 2n+1)$. 
Therefore $\Pi$ is a submersion at $j^r\Gamma(0)$. 
The type $\bbb$ of $\Pi(j^r\Gamma)$ at $t = 0$ for sufficiently large is given by 
$b_i = a_i - a_1 + 1$. Since the codimension $\Sigma_{a_1}$ in $J^r_{\rm{int}}(1, 2n+1)$ is given 
by $n(a_1 - 1)$, we have 
$$
{\mbox{\rm{Jet-codim}}}_{\mathcal D}(\aaa) = 
\sum_{i=1}^{n+1}(b_i - i) + n(a_1 - 1) 
= \sum_{i=1}^{n+1}(a_i - a_1 + 1 - i) + n(a_1 - 1) = \sum_{i=2}^{n+1}(a_i - i). 
$$
\QED

\

Third, similarly to above, we consider the jet space of ${\mathcal C}$-integral curves, 
$J^r_{\mathcal C}(I,  {\mathcal F}_{n+2}) \subset J^r(I, {\mathcal F}_{n+2})$ and 
$\Sigma_{\mathcal D}({\aaa})$ in 
$J^r_{\mathcal D}(I,  {\mathcal F}_{n+2})$. We denote its codimension by 
${\mbox{\rm{Jet-codim}}}_{\mathcal C}(\aaa)$. 

\bet
\label{codim-C}
The jet-codimension of the set of 
${\mathcal C}$-integral curves 
$\Gamma : I 
\to {\mathcal F}_{n+2}$ such that 
$\pi_1\circ\Gamma$ is of type $\aaa = (a_1, a_2, \dots, a_{n+1})$, is given by 
$$
{\mbox{\rm{Jet-codim}}}_{\mathcal C}(\aaa)  =  a_{n+1} - (n+1) \ = \ s(\aaa) - s(\aaa'), 
$$
where $\aaa' = (a_1, a_2, \dots, a_{n})$. 
\ent

\Proof 
As is explained in \S \ref{Pseudo-contact and canonical distributions.}, a ${\mathcal C}$-integral curve 
is described by the components $x_j^{j-1}(t), (1 \leq j \leq n+1)$. 
In fact, by projecting to these components, 
we see the diffeomorphism between the fiber $J^r_{\mathcal C}(I,  {\mathcal F}_{n+2})_{(t, f)}$ over 
$(t, f) \in I \times {\mathcal F}_{n+2}$ of the jet bundle 
and the ordinary jet space $J^r(1, n+1)$. 
To get the formula on ${\mbox{\rm{Jet-codim}}}_{\mathcal C}(\aaa)$, let 
$\Gamma(t) = (x_{ij}(t))$ be a ${\mathcal C}$-integral curve for the coordinates introduced in 
\S \ref{Pseudo-contact and canonical distributions.}  
through the origin at $t = 0$. 
Then we have, for the order at $t = 0$, 
$$
\ord x_{i}^{\ j} = \ord x_{i}^{\ j+1} + \ord x_{j+1}^{\ j}, \quad (0 \leq j, j+1 < i). 
$$
Therefore we have $\ord x_i^{\ 0} = \sum_{1\leq j \leq i} \ord x_j^{\ j-1}$ and hence 
$\ord x_i^{\ 0} - \ord x_{i-1}^{\ 0} = \ord x_i^{\ i-1} \geq 1, (2 \leq i \leq n+1)$. 
Then the type of $\pi_1\circ\Gamma$ is of type at $t = 0$ if and only if  
$\ord x_{i}^{\ 0} = a_i, (1 \leq i \leq n+1)$. The condition is equivalent to that 
$\ord x_i^{\ i-1} = a_i - a_{i-1}, (1 \leq i \leq n+1)$. 
Regarding the codimension in  $J^r(1, n+1)$, we have 
$$
{\mbox{\rm{Jet-codim}}}_{\mathcal C}(\aaa) = 
\sum_{i=1}^{n+1} (\ord  x_i^{\ i-1} - 1) = a_{n+1} - (n+1). 
$$
\QED

\ber
{\rm 
If $\gamma = (1, x_1^{\ 0}, \dots, x_{n+1}^{\ 0}) : I \to \R^{n+2} \setminus \{ 0\}$ 
at $t \in I$ is $\aaa = (a_1, a_2, \dots, a_{n+1})$, then the dual curve $\gamma^* = 
(1, x_{n+1}^{\ n}, x_{n+1}^{\ n-1}, \dots, x_{n+1}^1, x_{n+1}^{\ 0}) : I \to \R^{n+2} \setminus \{ 0\}$ is of type 
$$
\aaa^* = (a_{n+1} - a_n, a_{n+1} - a_{n-1}, \dots, a_{n+1} - a_1, a_{n+1}). 
$$
(Arnold-Scherbak's theorem \cite{O.P.Scherbak1}). 

}
\enr

As a consequence, we observe that 
$$
{\mbox{\rm{Jet-codim}}}_{\mathcal C}(\aaa) \leq {\mbox{\rm{Jet-codim}}}_{\mathcal D}(\aaa) 
\leq {\mbox{\rm{Jet-codim}}}(\aaa).  
$$
Since the transversality theorem (\cite{Mather}) hold, we see a curve of type $\aaa$ at a point 
appear generically if $s(\aaa) \leq 1$ in the class of curves in $P^{n+1} = \Gr(1, \R^{n+2})$, and 
$\aaa = (1, 2, \dots, n, n+1), (1, 2, \dots, n, n+2)$. 
Moreover, a curve of type $\aaa$ at a point 
appear momentarily in a generic one-parameter family of curves in $P^{n+1}$ 
if $s(\aaa) \leq 2$. The list is given by 
$$
(1, 2, \dots, n, n+1), (1, 2, \dots, n, n+2), (1, 2, \dots, n, n+3), (1, 2, \dots, n+1, n+2). 
$$

For adapted framed curves we have: 
\bet
\label{adapted}
For a generic one-parameter family of integral curves $\Gamma_\lambda : I \to {\mathcal F}_{n+2}$  
($\lambda \in J$, $J$ being a one-dimensional manifold) 
to the pseudo contact structure ${\mathcal D}$ 
on the flag manifold ${\mathcal F}_{n+2}$, the type of $\gamma_\lambda = \pi_1\circ \Gamma_\lambda$ 
and $\widehat{\gamma}_\lambda = \pi_{n+1}\circ\Gamma_\lambda$ 
at any point in $I$ for any parameter $\lambda \in J$ is one of the following list: 
$$
(1, 2, \dots, n, n+1), (1, 2, \dots, n, n+2), (1, 2, \dots, n, n+3), (1, 2, \dots, n+1, n+2), (2, 3, 4) (n = 2). 
$$
\ent

\Proof 
For ${\mathcal D}$-integral curves, the transversality theorem holds. In fact in \cite{Ishikawa6} 
the transversality theorem for integral curves to contact structure and the pseudo contact structure 
is the pull-back by the submersion $\pi : {\mathcal F}_{n+2} \to {\mathcal I}_{n+2}$ of the 
contact structure on the incident manifold ${\mathcal I}_{n+2}$. 
By Theorem \ref{codim-D}, we see ${\mbox{\rm{Jet-codim}}}_{\mathcal D}(\aaa) \leq 1$ 
if and only if $\aaa = (1, 2, \dots, n, n+1), (1, 2, \dots, n, n+2)$. 
Moreover ${\mbox{\rm{Jet-codim}}}_{\mathcal D}(\aaa) 
\leq 2$ if and only if $\aaa$ is one of the above list. Therefore we have the result. 
\QED

\

For osculating framed curves we have: 
\bet
\label{osculating}
For a generic one-parameter family of integral curves $\Gamma_\lambda : I \to {\mathcal F}_{n+2}$, 
($\lambda \in J$)   
to the canonical structure ${\mathcal C}$, 
the type of $\gamma_\lambda = \pi_1\circ \Gamma_\lambda$ 
and $\widehat{\gamma}_\lambda = \pi_{n+1}\circ\Gamma_\lambda$ 
at any point in $I$ for any parameter $\lambda \in J$ is one of the following list: 
$$
\begin{array}{c}
(1, 2, \dots, n, n+1), (1, 2, \dots, i, i+2, \dots, n+2) (0 \leq i \leq n), 
\vspace{0.2truecm}
\\
(1, 2, \dots, i, i+2, \dots, j, j+2, \dots, n+3) (0 \leq i < j \leq n+1).
\end{array}
$$
\ent
\Proof 
It suffices to note that the transversality theorem for 
integral curves $I \to {\mathcal I}_{n+2}$ for the contact structure (\cite{Ishikawa6}) implies 
the transversality theorem for ${\mathcal C}$-integral curves 
$I \to {\mathcal F}_{n+2}$. 
Then we have the required result by Theorem \ref{codim-C}. 
\QED

\section{Singularities of envelopes and their bifurcations.}
\label{Singularities of envelopes and their bifurcations.}

Let $\gamma : I \to X$ be a framed curve with framing $(e_1, \dots, e_{n+1})$. 
Then the {\it envelope} $E(\gamma)$ of $\gamma$,  
generated by the family of tangent hyperplanes $e_{n+1}^{\perp}(t) \subset T_{\gamma(t)}X$,  
is defined as follows (\cite{Thom}): Take the frame-dual $\widehat{\gamma} = e^{n+1} : I \to Y$ 
and take the fiber product 
$$
W := \{ (t, z) \in I \times Z \mid \widehat{\gamma}(t) = \pi_2(z) \} 
$$
of $\widehat{\gamma} : I \to Y$ and $\pi_2 : Z \to Y$. Then $W$ is an $(n+1)$-dimensional manifold. 
The envelope $E(\gamma)$ is defined as the set of critical values of the projection $\Pi_1 : W \to X$ defined 
by $\Pi_1(t, z) = \pi_1(z)$, for $(t, z) \in W$. 
The lifting $\widetilde{\gamma} : I \to G \subset \widetilde{\mathcal F}_{n+2}$ projects to 
the integral lifting $\overline{\gamma} : I \to Z \subset \widetilde{\mathcal I}_{n+2}$, 
$\overline{\gamma}(t) = (e_0(t), e_{n+1})(t)$, to ${\mathcal D}$. 
Note that $\pi_2\circ \overline{\gamma} = \widehat{\gamma}$. 
Moreover 
if we consider the osculating hyperplanes to $\widehat{\gamma} : I \to Y \subset \widetilde{\Gr}(1, \R^{n+2})$, (resp. $\widetilde{\Gr}(1, \R^{1, n+1})$), we get the dual curve 
$\widehat{\gamma}^* : I \to P^{n+2}$ (resp. $I \to P^{1, n+1}$), forgetting the orientation 
if necessary. In fact, $\widehat{\gamma}^* = \pi_1\circ\Gamma$ for the ${\mathcal C}$-integral lift 
$\Gamma : I \to {\mathcal F}_{n+2}^{0}$ of $\widehat{\gamma}$ constructed by associated osculating flags to 
$\widehat{\gamma}$, with respect to $\pi_{n+1} : {\mathcal F}_{n+2}^{0} \to P^{n+1}$ (resp. $P^{1, n}$). 

\

\noindent
{\it Proof of Theorem \ref{envelope}.} 
By Theorem \ref{adapted}, the type of 
a curve $\widehat{\gamma}_\lambda : I \to Y$ in a generic one-parameter family of framed curves 
at a point $(t, \lambda)$ is one of ${\rm (I)}: (1, 2, 3), {\rm (II)}: (1, 2, 4), {\rm (III)} : (1, 3, 4), {\rm (IV)} : (1, 2, 5)$ or ${\rm (V)} : (2, 3, 4)$. 

Then we have the normal forms of singularities by the classification results in \cite{Ishikawa1}. 
\QED

\

Moreover, by Theorem \ref{osculating}, we have immediately: 
\bet
For a generic one-parameter family $\gamma_\lambda$ of osculating framed curves in $E^3, S^3, H^3$, 
the frame dual $\widehat{\gamma}_\lambda$ 
at a point $(t, \lambda)$ has one of type in the list 
$$
(1, 2, 3); (1, 2, 4), (1, 3, 4), (2, 3, 4); (1, 2, 5), (1, 3, 5), (1, 4, 5), (2, 3, 5), (2, 4, 5), (3, 4, 5). 
$$
Corresponding to each type in the above list, 
the dual curve $\widehat{\gamma}_\lambda^*$ turns to be of type 
$$
(1, 2, 3); (2, 3, 4), (1, 3, 4), (1, 2, 4); (3, 4, 5), (2, 4, 5), (1, 4, 5), (2, 3, 5), (1, 3, 5), (1, 2, 5).
$$
\ent

By \cite{Ishikawa1}, the diffeomorphism class of the envelope 
$E(\gamma_\lambda)$ is determined in the case
$$
(1, 2, 3), (2, 3, 4), (1, 3, 4), (1, 2, 4), (3, 4, 5), (1, 2, 5). 
$$ 
Moreover in any case the topological class of $E(\gamma_\lambda)$ is determined (\cite{Ishikawa2}). 
We will describe in the forthcoming paper in detail, 
the topological bifurcations of envelopes for osculating framed curves of type 
${\mbox{\rm{Jet-codim}}}_{\mathcal C}(\aaa) \leq 2$, namely, for 
$$
(2, 4, 5), (1, 4, 5), (2, 3, 5), (1, 3, 5), (1, 2, 5). 
$$ 

\ber
{\rm
Though the list of singularities is common for all of three geometries, the geometric characters are of course 
distinguished. For instance, in the case $n = 2$ and for an adapted frame $e_0' = e_1$, 
we have the structure equation, under the arc-length derivative, 
$$
\left\{
\begin{array}{cccccc}
e_1' & = & - \delta e_0 & & +  \kappa_1 e_2  & + \kappa_2 e_3,  \\
e_2' & = & & - \kappa_1 e_1  &  & + \kappa_3 e_3,  \\
e_3' & = & & - \kappa_2 e_1 & - \kappa_3 e_2, & 
\end{array}
\right. 
$$
where $\delta = 0, 1, -1$ for $X = E^3, S^3, H^3$ respectively (\cite{IL}). 
In general, we characterise the type $\aaa$ of the frame-dual $\widehat{\gamma}$ for a framed curve $\gamma$ 
by polynomials, distinguished in each gemetry, 
of geometric invariants of $\gamma$ and their derivatives up to order $a_{n+1} - 1$. 
For $E^3$, see \cite{Ishikawa6}. 
}
\enr

{

}

\

\noindent Goo ISHIKAWA,  \\ 
Department of Mathematics, 
Hokkaido University, 
Sapporo 060-0810, JAPAN. 
\\
{\vspace{-0.7truecm}}
\begin{verbatim}
E-mail : ishikawa@math.sci.hokudai.ac.jp
\end{verbatim}

\end{document}